# Recency Bias in the Era of Big Data: The Need to Strenthen the Status of History of Mathematics In Nigerian Schools

**Abah, J.A.**
Department of Science Education
University of Agriculture
P.M.B. 2373 Makurdi
Benue State, Nigeria.
abahjoshua@uam.edu.ng

## ABSTRACT

The amount of information available to the mathematics teacher is so enormous that the selection of desirable content is gradually becoming a huge task in itself. With respect to the inclusion of elements of history of mathematics in mathematics instruction, the era of Big Data introduces a high likelihood of Recency Bias, a hitherto unconnected challenge for stakeholders in mathematics education. This tendency to choose recent information at the expense of relevant older, composite, historical facts stands to defeat the aims and objectives of the epistemological and cultural approach to mathematics instructional delivery. This study is a didactic discourse with focus on this threat to the history and pedagogy of mathematics, particularly as it affects mathematics education in Nigeria. The implications for mathematics curriculum developers, teacher-training programmes, teacher lesson preparation, and publication of mathematics instructional materials were also deeply considered.

## 1. INTRODUCTION

Mathematics is an embodiment of tradition. As a field of study, mathematics is the conglomeration of centuries of efforts structurally pieced together from diverse facet of human endeavor and from across different cultures around the globe. This structure of knowledge is communicated over the ages to students in schools as a composite subject, with its practice and evolution deeply embedded in all successful (and sometimes unsuccessful) human advancement. The fact of the relevance of mathematics to life elevates the status of the mathematics educator to paramount. In expressing this truth, Pepin (1999) rightly asserts that mathematics educators are charged with explaining the processes of growth of mathematical knowledge, monitoring the processes of mathematical discovery in the making of the mathematician, and provoking such processes in teaching. Mathematics teachers are therefore expected to be the custodians of the rich history and traditions of the discipline.

Mathematical creativity ensures the growth of the field of mathematics as a whole (Sriraman, 2004). The ingenuity of mathematicians past are today manifesting in the indispensability of modern science and technology. Their painstaking strides have woven a rich mosaic of conceptions of the nature of mathematics, ranging axiomatic structures to generalized heuristics for solving problems (Dossey, 1992). In its earliest days (spanning the era of the Egyptian and Greek thinkers) mathematics was often bound up with practical inquiry (Krantz, 2007). However, as time progresses, it became more balanced and more realistic to consider that mathematics exists as a body of truths about relationships between abstract entities and structures (Hatcher, 2008).

These abstract relationships are reflected or instantiated, in various forms and at different levels, in the concrete structure of the physical world. Even in its most advanced form, the subject is rooted in reality with applications in everyday life. Mathematics is vital to business, agriculture, medicine, psychology and several fields of human endeavour.

The permeability of mathematics made emphasis on the pedagogy of the subject a central concern to educators. The burden of the sustainability of mathematical knowledge rests squarely on the shoulders of mathematics teachers, whose methodology of communicating facts determine both immediate and future outcomes for students. It is obvious that mathematics in some sense has a common language: a language of symbols, technical definitions, computations, and logic (Thurston, 1994). The rectitude of mathematics requires students to articulate sound mathematical explanations and always justify their solutions. Currently, focus is shifting from teacher-centred instructional strategies that stunt the growth of mathematical mastery to student-driven mathematical intuition, mathematical thinking and real-world problem solving (Legner, 2013). Linking the teaching of mathematical concepts to the histories behind them has been associated with setting students on the path of mathematical discovery, irrespective of the method of instruction (Yevdokimov, 2006).

Over recent decades, advances in digital technology have provided several forms of augmentation to the traditional mathematics classroom. Today, powerful information and communications technology (ICT) tools and platforms at the disposal of the mathematics teacher are helping transform education in terms of its contents, methods and outcomes (Iji & Abah, 2016). Both teachers and students of mathematics now have access to vast amount of resources at their fingertips wherever they are. Indeed there have been a massive revolution in the way mathematics learners study and do research. This digital revolution is being fueled by the increasing broadband penetration and ubiquity of smartphones among the new genre of learners who are successively imbibing a trendy culture of leisure and school work. Tons of digital content are being churned out by different outlets daily over the internet, with the Computer Science Corp (2012) forecasting an annual increase of 4300% by the year 2020 when the information super-highway will be holding an estimated milestone traffic of 35 Zettabyte of data. This phenomenon of data explosion playing out in the current era is termed "Big Data".

Big data has been described as the incomprehensibly large worlds of information which have been on the rise due to explosion in mobile networks, cloud computing and new technologies (Bollier, 2010). According to Boyd and Crawford (2012), Big Data is less about data that is big than it is about a capacity to search, aggregate, and cross-reference large data sets. The interconnection of databases and development of powerful new software tools for inference making are resulting in a new kind of knowledge infrastructure freely available for patronage by every inquisitive mind. Educational establishments are active players in this knowledge economy. But like the Aspen Institutes Executive Director Charles M. Firestone asked in his forward to the Bollier (2010) report, "does the ability to analyze massive amounts of data change the nature of scientific methodology? Does Big Data represent an evolution of knowledge, or is more actually less when it comes to information on such scales?" (p.viii). The answers to these questions are fast constituting a heated debate among scholars, and obviously the promises outweighing the perils.

When eyeballed in the context of mathematics education, the existence of Big Data is immensely widening the scope of coverage for both teachers and students. It can be argued that mathematicians have more access to inexhaustible resources now than ever before. The attainable horizon of research in the discipline is now unendingly scalable. However, with respect to the study of history of mathematics, the issues become more confounding and complex than outright. With more new information being added daily to educational resources, focus on integrating elements of history into normal mathematics classroom activities may suffer some setbacks, not excluding the tendency to resort to more recent history. If this bias becomes pervasive in mathematics education practice, the very objectives of the inclusion of elements of history will be gravely endangered.

The psychological tendency to reference more recent information at the expense of older, detailed and more relevant data, rightly termed "Recency Bias", may be gradually eroding the foundations of mathematical knowledge. The chain of events weaved up in the history of basic mathematical concepts are fundamental to enriching students understanding of school mathematics. It is against this backdrop that this work seeks to lay out the issues contending with the handling of history of mathematics, particularly in Nigeria. The study intends to establish the current status of history of mathematics in Nigeria and the need for various stakeholders in mathematics education to rise to the challenge of saving this aesthetic component of mathematics from an impending demise.

## 2. RELEVANCE OF HISTORY OF MATHEMATICS

Mathematics is proud of its history, confident in its traditions and certain that the truth it pronounced are real truths applying to the real world (Wells, 2006). There is nothing as dreamy and poetic, nothing as radical, subversive, and psychedelic, as mathematics, since it is the purest of the arts, as well as the most misunderstood (Lockhart, 2002). Mathematics, as taught in school has reduced the rich and fascinating adventure of the imagination to a sterile set of facts to be memorized and procedures to be followed.

History of mathematics adds aesthetic value to the teaching and learning of mathematics. Teaching mathematics without adequately incorporating the elements of its history leads to removing the creative process and leaving only the results of that process. This will automatically lessen the chance of any real engagement with the subject. History of mathematics places mathematics in real time and location, and blending it into classroom instruction, according to Lockhart (2002), enables students to pose their own problem, make their own conjectures and discoveries, to be wrong, to be creatively frustrated, to have an inspiration and to cobble together their own explanations and proofs. These rich experiences are what makes mathematicians, both ancient and modern, admirable characters in the society. Denying students the opportunity to engage in these activities implies denying them mathematics altogether.

Mathematics should be thought of as a living organism with its peculiar long history and a vivid present. Considering this, mathematics is expected to be taught with reference to a historical, epistemological and cultural approach. Vallhonesta, Esteve, Casanova, Puig-Pla and Roca-Rosell (2015) stated that history of mathematics can be used implicitly and explicitly to enrich mathematics education. On the one hand, the history of mathematics can be used as an implicit resource in the design of activities to adapt some standard concepts, and to prepare problems and auxiliary sources. On the other hand the history of mathematics can be used in an explicit way to direct and propose research works using historical materials to aid students' understanding of mathematical concepts. History of mathematics helps to develop the analytic and synthetic thought process of students. The ideas generated from historical attachments to mathematical concepts can truly create conflictive situation in which students are encouraged to reflect upon the rules that define their actions when dealing with the concepts (Bernardes & Roque, 2015). This approach, results in both sensitive and historical thinking, mediated by bodies, signs, artefacts and cultural meanings (Guillemette, 2015). The approach gives rise to a non-mentalist conception of thought.

Generally, students' reflections about the nature of mathematics through history evidently results in enhanced mathematical literacy, increased psychological motivation, and linguistic and transverse competencies. The usage of original sources has been shown to be effective in the teaching and learning mathematics (Guillemette, 2015). In view of the usefulness of history of mathematics, content selection must be a careful effort on the part of the teacher in order to achieve an unbroken connection of related facts.

## 3. RECENCY BIAS ARISING FROM BIG DATA

Recency, basically, is a cognitive bias that convinces one that new information, which is more recent, is more valuable and important than older information. According to Townsend (2014), this may be true, but it is not necessarily true. Kanasky (2014) attributes recency effect to a simple enhancement of short-term memory due to recent exposure to information. Sousa (2006) related that Ebbinghaus published the first studies on this phenomenon in the 1880s, and as such, recency bias itself is not a new discovery.

Cases of recency bias are well established in the fields of Educational Psychology, Law, and Behavioural Finance. Jones and Sieck (2003) having considered advantageous and disadvantageous scenarios of recency effect, suggested that regardless of the mechanisms responsible, it is instructive to consider their role in relationship to normative and descriptive behaviour. Barbosa, de Lima-Neto, Evsukoff and Menezes (2015) found that when it comes to visitation patterns, humans are extremely regular and predictable, where recurrent travels responsible for most of the movements. Teachers usually circumvent the recency effect in the classroom often by applying review, warm-ups, recap of previous knowledge, practice, and summary. Lawyers, for instance, call their strongest witness last to make a lasting impression. Advertisers also apply the recency factor by ensuring the last portion of their promotion creates the desire to purchase their product.

Recency bias is often classified as one of many different decision-making errors and closely related to the human tendency towards cognitive ease as against the hard work of thinking, analyzing and deciding. It is the habitual laziness factor in recency bias that distinguishes it from the seemingly similar recency effect which is dictated by limitation (SKYbrary, 2016). This fact links recency bias to choice, which is in turn heavily dependent on the amount of available possibilities. Considering how limited abilities to process information affect choice behaviour (Salant, 2011), the current era of Big Data certainly stretches the range of available possibilities to an unimaginable extent. While it is advisable to weigh as much facts and contextual information as possible before making decisions, is not always practicable in the decision making process of lesson preparation and individual study. Those in charge of selecting  instructional aid on both online and offline media are subject to recency bias and history of mathematics as an element of mathematics instruction tends to be adversely affected considering the ease of neglecting pre-modal historical data by the non-professional eye on the grounds of perceived irrelevance or age.

Processes of globalization resulting in Big Data have produced a strong cross-cultural perspective on how people recall and historical narratives (Gonzalez-Castro, 2006). The chain of events which constitute the body of communicated history is required to be transmitted systematically. However, in the case of history of mathematics, teachers' departure from relating specific origins of concepts and attempts by ancient mathematicians in providing solutions to mathematical problems is contributing to the mystification of the subject. Students are erroneously led to believe that mathematics is highly abstract because concepts treated in the classroom are left hanging without meaningful beginning and representation in real life. An illustration can be cited of the case of analytic geometry in which most instructional contents dwell only on deriving and resolving equations without at least going backward to the significant contributions of Rene Descartes, or even backtracking further to early attempts by ancient Egyptian and Greek mathematicians.

Recency bias in the development of mathematical contents occasioned by the multiplicity of sources in this era of Big Data is making mathematics more mechanic than realistic. Mathematics is mechanical if it is devoid of history and applications. Teachers of mathematics, being aware of this fact, are expected to be seasoned editors of materials with interest in promoting the aesthetic value of the subject by transcending cultures and ages to connect the contents they are required to deliver.

## 4. THE NEED TO BEGIN AT THE BEGINNING IN NIGERIA

The joy in mathematics can be felt in discovering new mathematics, rediscovering old mathematics, learning a way of thinking from a person or text, and finding a new way to view an old mathematical structure (Thurston, 1994). The teacher of mathematics, irrespective of the level of education, is at the centre of communicating a sound mathematical experience to students. It is the duty of such teacher to design mathematics instruction that will dynamically utilize aspects of concept history to enrich the learning experience.Although, Nigeria is on the African continent known for its rich heritage of mathematics culture, indigenous studies into the history of its mathematics are scanty and still hard to come by. However, pockets of studies into the integration of cultural artefacts in the study of mathematics affirm the relevance of history of mathematics in the communication of mathematics in Nigeria (Shuaibu, 2014). Full empirical deployment of history of mathematics in the classroom is still unavailable from Nigeria. There is the indication that little elements of history are being embedded in classroom instruction in mathematics considering how General History fared in the development of curriculum in Nigeria.

A glimpse through the literature available for the teaching and learning of mathematics attest to the poor handling of history of mathematics in Nigeria. Techers who did not wake up to the discovery of the relevance of history of mathematics are left groping with the teaching of the subject only as they, in turn, were taught. In the absence of innovation, students are mechanically put through the age-long myth of difficulty of mathematics by grinding through sets of solutions without any explanation of the mechanisms involved. Heightened mathematics phobia arising out of this status quo drove away a good number of students from pursuing mathematics related careers. The few who could attain mathematical proficiency out of personal effort are viewed as privileged and gifted, making it appear like it is unnatural to understand mathematics. There is obviously a serious need for mathematics educators to begin at the beginning. Teachers should bring in historical connections that give students the opportunity to model mathematical problem in its existential context, to carry out analysis of the learning materials and discover mathematical properties that are completely new (Yevdokimov, 2006). One of the means of demystifying mathematics is for teachers to introduce original sources to their classroom discussion. The archetypal texts of early attempts of ancient mathematicians can be used as starting points for classroom discussion. A prominent example of original source instructional materials is the Chinese text: The Nine Chapters in Mathematical Arts, which has ties to most mathematical procedures of the current era. Systems of linear equations, trigonometry and geometry are aspects of school mathematics that could be spiced up by referencing this ancient text as a starting point.

As soon as students attain the formal operational stage (11-16 years), they are ready to be confronted with historical materials that foster hypothetical deductive reasoning. The support of carefully articulated elements of history can spark up reflective thought process from concrete considerations, building up to more abstract reasoning (Joubish & Khurram, 2011). Such articulation of materials must be free of prejudices arising from Big Data.

## 5. IMPLICATIONS FOR MATHEMATICS EDUCATION IN NIGERIA

There is accumulating evidence that recent prior knowledge about expectations plays an important role in perception (Raviv, Alhissar & Loewenstein, 2012). As the Internet continues to be interwined with peoples' daily lifestyle, it is negatively shaping the way information is being processed and interpreted. Weyers (n.d) observed that the way teachers and students are currently using the Internet is reducing the desire to be inquisitive, think, comprehend, and ultimately retain information. Specifically, the era of Big Data is increasingly contributing to the display of recency bias with respect to the handling of history of mathematics. The implications of this didactic discourse are considered here from some key perspectives. These are as it relate to curricular emphasis, teacher preparation, reformation of publications of basic mathematics literature, and teacher training programmes. Literature suggests that most creative individuals tend to be attracted to complexity often arising out of past problems, an existence of which most school curricula has very little to offer (Sriraman, 2004).

Mathematics curricular efforts rarely throw open problems with the sort of underlying mathematical structure that would warrant students having a prolonged period of engagement and the independence to formulate solutions. Mathematics curricular efforts in Nigeria should seek to attain the higher aim of mathematics education, which is to establish in students clarity of thought and the ability to pursue assumptions to logical conclusions. Re-integrating elements of mathematical history to the background of instructional content could be a safe starting point for both developers and implementers of the curriculum.

The Nigerian Education Research and Development Council (NERDC) which is charged with this onerous regulatory function must be repositioned to effect sound instructional delivery in mathematics across the country. Developers must stand up against the idea of recycling mathematics curriculum with no historical perspective and should put in place mechanisms for restoring discovery and exploration, natural and meaningful problem solving, and factual simplicity. Practical cognitive artefacts like early classical texts and historical approaches to solving mathematical problems can be recommended by the curriculum. If these are passed down with the publication of textbooks and other instructional materials, curious students can lay hold on their benefits to drive more personalized learning.

Mathematics is a slow, contemplative process which takes time to produce and a skilled teacher to even recognize one. Lockhart (2002) posit that it is unacceptable to have mathematics teachers who know nothing of the history and philosophy of the subject, nothing about past and recent developments, nothing in fact beyond what they are expected to present to their unfortunate students. Only a high level of co-ordination on the part of the mathematics teacher can prevent mathematics instruction from being reduced to mere data transmission, a process in which there is no sharing of excitement and wonder. The teacher of mathematics is expected to first develop himself at the personal level in the art of building openness and honest intellectual relationship with students by leveraging on the connecting power of rudimentary mathematical history in instruction delivery. Every mathematics class should be loaded with rich learning experiences based on the teachers' ability to connect the past to the present and recognize information without prejudice to distance from the current era. In handling the distractions of Big Data, the mathematics teacher, like a professional driver, must avoid the temptation of being overwhelmed by the speed of visual data hitting him, by adopting expert vigilance and focusing on a far-away point which is the goal of his duty (Askitas, 2016). Without confident, knowledgeable teachers in the classroom, improvement in mathematics education is certainly impossible (Andrews, 2012).

The foregoing implication partially re-directs responsibility to teacher-training programmes in Nigeria. Universities and colleges of Education mounting programmes in mathematics education must emphasize the relevance of history of mathematics in the pedagogy of mathematics. Professionals at this level should inculcate in both in-service and pre-service teachers the notion that elements of history of mathematics constitute effective tools to help in the understanding of mathematics as useful, dynamic, humane, interdisciplinary and heuristic science (Vallhonesta et al, 2015). Apart from the teaching of the usual mandatory two-credit-units course – History of Mathematics, mathematics educators should practically engraft history into their own classroom discussions, thereby setting a good example for the teachers-in-training to follow. Nigerian mathematics education practitioners should set the pace in developing instructional strategies that make use of original sources from across the globe. In already developed forms, these tools can easily be picked up by mathematics teachers at the secondary and primary education levels for onward blended utilization. The widespread dissemination of mathematical knowledge is expected to present a more humanistic idea of what mathematics is by reporting mathematical history without breaking its genealogy. Basic mathematics literature available for use at lower educational levels in Nigeria is expected to position mathematics in time by adequately linking concepts to their origins. Publishers of textbooks, the National Mathematical Centre (NMC), Abuja, and the Mathematical Association of Nigeria must break away from the practice of producing literature that enforces a rigid view of mathematics. Online publications intended for classroom consumption should also cover essential aspects of historical linkages to key mathematical concepts. Materials in study books should be presented first from the concrete level, before gradually increasing in abstraction.

## 6. CONCLUSION

This discourse has attempted to stress the danger Big Data poses to the inclusion of wholesome elements of history in mathematics instructional delivery. In the presence of so much information, even mathematics teachers as guides in the teaching and learning process can succumb to the tendency to consider more recent history at the expense of a complete picture of origin of concepts. The need to reposition history of mathematics as an indispensable component of mathematics instruction was also highlighted. The Nigerian situation was given peculiar attention to emphasize the role history of mathematics can play in restoring the aesthetic value of mathematics among all levels of students in the country. The implications of the realities unraveled by this study were pointedly considered along its significance to Nigerian teacher-training institutions, mathematics teachers and students, curriculum developers and content publishing outlets. Like other herald of this timely tintinnabulation, this study has brought to light the need to preserve academic creativity and critical thinking. Indeed there is the need for selective modesty irrespective of the amount and speed of processing of available information. All stakeholders in mathematics education must pay wide attention to the chance of losing originality to recency bias as a result of Big Data.